\documentstyle[11pt]{article}
\begin{document}
\begin{center}
{\large \bf NUMERICAL PERFORMANCE OF 
ABS CODES 
\vskip2mm FOR NONLINEAR SYSTEMS OF EQUATIONS}
\footnote{Work partly supported by grant GA CR 201/00/0080
and by MURST  1999 Programma
di Cofinanziamento}

\end{center}

\begin{center}
{\bf E. Bodon
\footnote{Department of Mathematics, University of Bergamo, Bergamo 24129,
Italy (bodon@unibg.it)},
A. Del Popolo
\footnote{Department of Mathematics, University of Bergamo, Bergamo 24129,
Italy (delpopolo@unibg.it)},
L. Luk\v{s}an
\footnote{Institute of Computer Science, Academy of Sciences of the Czech
Republic,
Pod vod\'arenskou v\v e\v z\'\i\ 2, 182 07 Prague 8, Czech Republic
(luksan@uivt.cas.cz)} and
E. Spedicato
\footnote{Department of Mathematics, University of Bergamo, Bergamo 24129,
Italy (emilio@unibg.it)}}

\end{center}
\vskip5mm
{\bf 1. Introduction}
\vskip3mm
The nonlinear ABS algorithms solve a nonlinear system of algebraic
equations $F(x)=0$ where $F$ is a function $F:R^n\rightarrow R^n$ in
component form $f_1(x)=0, f_2(x)=0,..,f_n(x)=0, \: x\in R^n$.
The ABS algorithms can be viewed as a modification of  
the Newton method. Newon method constructs an iterative sequence 
of points $x_i, \: i=1,2,...$ by the relation $x_{i+1}=x_i+t_i$, where $t_i$
is the solution of the linear system $J(x_i)t_i=-F(x_i)$ or in component
form $a^T_k(x_i)t_i=-f_k(x_i), \: k=1,..,n$, where $J(x_i)$ indicates the
Jacobian matrix of $F$ evaluated at $x_i$, while $a_k(x_i)$ denotes the $k$-th
row of $J$ evaluated at $x_i$. The set of linear equations
$a^T_k(x_i)t_i=-f_k(x_i), \: k=1,..,n$ is solved by the ABS method for
linear systems after modifying the definition of $J$ and $F$. The motivation
is based on the fact that at the $k$-th step the ABS algorithms generate
the next approximation to the solution $t_i$ using the $k$-th equation and 
the vector $t_{i_k}$ which is a solution of the system formed by the first
$k-1$ equations. Since the new approximate point $x_i+t_{i_k}$ is expected to 
be closer to the solution $x^\star$ of $F(x)=0$, it can be used to evaluate
the $k$-th row of $J$ and the $k$-th component of $F$. The procedure is the
following.
\vskip3mm
Assing an initial vector $x_0$, set $m=1$, \\
($\star$) set $y_1=x_{m-1}$, $H_1=I$, for $k=1$ to $m$ compute: \\
$p_k=H_k^Tz_k$ where $z_k\in R^n$ is a parameter which gives the type of
ABS method,\\
$\beta_k=f_k(y_k)/\delta_k, \: \delta_k=p_k^Ta_k(y_k)$,
$y_{k+1}=y_k-\beta_kp_k$, $H_{k+1}=H_k-H_ka_k(y_k)p^T_k/\delta_k$, \\
let $x_m=y_{n+1}$,\\
if the stopping criterion is satisfied, stop; otherwise increment the
index $m$ by $1$ and return to ($\star$).
\vskip5mm
{\bf 2. The implemented codes}
\vskip3mm
The ABS codes that we have implemented based on the above algorithm 
use three choices of the available parameter $z_k$ that define
three well-known ABS methods for the linear case.

The first code, routine nl-huang1, and the second code, routine
nl-huang2, use the Huang method defined by $z_k=a_k(y_k)$, and the
modified Huang method defined by $z_k=H^T_ka_k(y_k)$ and $ \delta_k=p_k^Tp_k$.

In the routine nl-huang1 the matrix $H_k$ is updated by formula
$H_{k+1}=H_k-H_ka_k(y_k)p^T_k/\delta_k$. 
In one major iteration (consisting of a sequence of $n$ minor iterates such
that system $J(x)t_i=-F(x)$ is solved) the number of moltiplications is
$3/2 n^3$ for the Huang algorithm, $5/2 n^3$ for the modified Huang 
algorithm  and symmetry of $H_k$ is forced by computing only the elements on 
a triangular part of the matrix. The memory requirement is $n^2/2$.

In the routine nl-huang2 the vector $p_k$ is computed by the equivalent
formula $p_k=a_k(y_k)-P_{k-1}D^{-1}_{k-1}P^T_{k-1}a_k(y_k)$ in Huang method 
and $p_k=s_k-P_{k-1}D^{-1}_{k-1}P^T_{k-1}s_k$, 
$s_k=a_k(y_k)-P_{k-1}D^{-1}_{k-1}P^T_{k-1}a_k(y_k)$ in the modified Huang 
method, where $P_{k-1}=(p_1,..,p_{k-1})$, 
$D_{k-1}=(\delta_1,..,\delta_{k-1})$.
In one major iteration the number of moltiplications is $n^3$ for the Huang 
algorithm, $2 n^3$ for the modified Huang algorithm, the memory is $n^2$. 
This code requires less operations but more memory than nl-huang1.

In the third code, routine nl-ilu, the parameter choices are the same of 
the implicit LU algorithm with column pivoting, that is $z_k=e_j$ where $j$
is such that $ |e_j^TH_ka_k(y_k)|=max\{ |e_l^TH_ka_k(y_k)|, l=1,..,n\}$.
A major iteration requires $n^3/3$ multiplications and $n^2/2$ memory
(memory requirement can be reduced to $n^2/4$ by using a different 
implementation). 

In a major iteration the number of component function evaluations is $n$ 
and the number of element jacobian evaluation is $n^2$ for all the codes.

The jacobian is singular if at least one row $k$ is linearly dependent, in
which case in all the algorithms the iteration continues setting 
$y_{k+1}=y_k$. In the theoretical ABS method the test for linearly dependence
is $H_ka_k=0$, in the present 
codes the test is done on the value of $\delta_k$ and
the zero is substituted by a tolerance $t$ or $t*\|a_k\|$ where $t$ is given
by the user.

The programs have four stopping criteria. The first criterion depends on a 
measure of the residual $\|F(x_i)\|_\infty \le eps$, the second criterion is 
based on a measure of the relative distance between consecutive iterates
$\|x_{i+1}-x_i\|_\infty \le tol\|x_i\|_\infty$, the third criterion
terminates the iterations if no progress occurs in decreasing the function
after a fixed number $n_s$ of steps. The last criterion terminates the 
iterations after a fixed number $itmax$ of steps. The parameters
$t, \: eps, \: tol, \: n_s, \: itmax$ are given by the user.

Optionally, if the norm of $F$ is increasing at an iteration $i$, the 
implemented algorithms can 
use a line search technique. Setting $\bar{x}=x_i$, the 
line search technique reduces the width of the interval between the points 
$x_{i-1}$ and $\bar{x}$ by $1/2$ (the new point is 
$\bar{x}=(x_{i-1}+\bar{x})/2$) and so on, until a point $\bar{x}$ is 
found where the norm of $F(\bar{x})$ decreases, or a maximum number 
$n_{half}$ of halvings is achieved, where $n_{half}$ is given by the user.

The codes are written in Fortran 77 in single and in double precision.
\vskip5mm
{\bf 3. The numerical experiments}
\vskip3mm
Some numerical results are given in the following 
tables. The symbols mod.huang1, mod.huang2, 
implicit lu, refer to the codes nl-huang1, nl-huang2, nl-ilu, the symbols  
m.hua1 line search, m.hua2 line search, imp.lu line search, mean that the 
line search technique was applied in the algorithms nl-huang1, nl-huang2, 
nl-ilu, where the function norm grows at some iterations. The symbol
$\|F\|_\infty$ means the minimum value found by the algorithm,
$it_{ \|F\|_\infty}$ the iteration where $\|F\|_\infty$ was found,
$it$ the number of performed iterations, $time$ the execution time.
The symbols in the brackets have the following meanings: (x) the second
stopping criterion is verified, (div) divergence, stop because $\|F\|_\infty$
increases during a fixed number of steps, (o) oscillations, stop because 
there is no improvement in the last $n_s$ steps. Lack of marks means 
that the first convergence test is satisfied.

The experiments are obtained setting the tolerance parameters $t=1.e-6$,
$eps=1.e-6$, $tol=1.e-10$ for the single precision case and $t=1.e-15$,
$eps=1.e-15$, $tol=1.e-18$ for the double precision case.

\smallskip 
We have used the following test functions:
\begin{itemize}
\item Rosenbrock function ($n=2$ and the extended Rosenbrock function, where, for 
\vskip2mm
$i = 1,2.... n/2$\\
$f_{2i-1}=1-x_{2i-1}$ \\
$f_{2i}=10(x_{2i}-x_{2i-1}^2)$ \\
$x_{0_{2i-1}}=-1.2, \quad x_{0_{2i}}=1.0 $  \\   
\item Powell singular function with $n=4$, where \\
\vskip2mm
$f_1=x_1+10x_2$\\
$f_2=5^{1/2}(x_3-x_4)$\\
$f_3=(x_2-2x_3)^2$ \\
$f_4=10^{1/2}(x_1-x_4)^2$ \\
$x_0=(3,-1,0,1)$\\
\item Brown almost linear function \\
\vskip2mm
$f_i=x_i+\sum_{j=1}^nx_j-(n+1)$ \\
$f_n=(\prod_{j=1}^nx_j)-1$ \\
$x_{0_i}=1/2, \quad i=1,2,..,n$ \\   
\item Schubert Broyden function \\
\vskip2mm
$f_1=(3-x_1)x_1+1-2x_2$ \\
$f_i=(3-x_i)x_i+1-x_{i-1}-2x_{i+1} \quad i=2,..,n-1 $ \\
$f_n=(3-x_n)x_n+1-x_{n-1}$ \\
$x_{0_i}=-1, \quad i=1,2,..,n   $  \\   
\end{itemize}
\vskip5mm
{\bf 4. References}
\vskip3mm
\newpage 

{\scriptsize
\begin{tabular}{llccccc}  
\multicolumn{2}{c}{ }&\multicolumn{5}{c}{single precision } \\  \\
function & method &$\|F\|_\infty$& &$it_{ \|F\|_\infty}$&$it$&$time$\\  \\
Brown almost linear n=4 & mod.huang1  &0.60e-6 & &3 &3 & 0.0 \\
  $x_0$                 & mod.huang2  &0.60e-6 & &3 &3 & 0.0 \\ 
                        & implicit lu &0.24e-6 & &5 &5 & 0.0 \\  \\
Brown almost linear n=4 & mod.huang1  &0.12e-6 & &3 &3 & 0.0 \\
  $1.1x_0$              & mod.huang2  &0.60e-7 & &3 &3 & 0.0 \\ 
                        & implicit lu &0.12e-6 & &5 &5 & 0.0 \\  \\
Brown almost linear n=4 & mod.huang1  &0.60e-6 & &4 &4 & 0.0 \\
  $10x_0$               & mod.huang2  &0.36e-6 & &4 &4 & 0.0 \\ 
                        & implicit lu &0.10e+1 &(x) &1 &2 & 0.0 \\  \\
Brown almost linear n=4 & mod.huang1  &0.42e-6 & &14 &14 & 0.0 \\
    $100x_0$     & m.hua1 line search &0.60e-6 & &9  &9  & 0.0 \\  
                        & mod.huang2  &0.36e-6 & &14 &14 & 0.0 \\ 
                 & m.hua2 line search &0.24e-6 & &9  &9  & 0.0 \\ 
                        & implicit lu &0.0 & &17 &17 & 0.0 \\  \\
Brown almost linear n=20& mod.huang1  &0.24e-6 & &41 &41 & 1.59 \\
    $x_0$        & m.hua1 line search &0.72e-6 & &7  &7  & 0.27 \\  
                        & mod.huang2  &0.95e-6 & &14 &14 & 0.61 \\ 
                 & m.hua2 line search &0.72e-6 & &9  &9  & 0.39 \\ 
                        & implicit lu &0.60e-6 & &14 &14 & 0.16 \\  
                 & imp.lu line search &0.72e-6 & &12 &12 & 0.16 \\ \\ 
Brown almost linear n=20& mod.huang1  &0.24e-5 &(o) &12 &15 & 0.55 \\
    $1.1x_0$     & m.hua1 line search &0.36e-6 & &9  &9  & 0.39 \\  
                        & mod.huang2  &0.25e-5 &(o) &7 &10 & 0.44 \\ 
                 & m.hua2 line search &0.48e-6 & &11  &11  & 0.49 \\ 
                        & implicit lu &0.25e-5 &(o) &6 &10 & 0.11 \\  
                 & imp.lu line search &0.36e-6 & &9 &9 & 0.11 \\ \\ 
Schubert Broyden n=10   & mod.huang1  &0.42e-6 & &4 &4 & 0.0 \\
  $x_0$                 & mod.huang2  &0.42e-6 & &4 &4 & 0.0 \\ 
                        & implicit lu &0.89e-6 & &3 &3 & 0.0 \\  \\
Schubert Broyden n=10   & mod.huang1  &0.48e-6 & &8 &8 & 0.5 \\
  $10x_0$               & mod.huang2  &0.48e-6 & &8 &8 & 0.5 \\ 
                        & implicit lu &0.54e-6 & &7 &7 & 0.5 \\  \\
Schubert Broyden n=50   & mod.huang1  &0.24e-6 & &4 &4 & 2.03 \\
  $x_0$                 & mod.huang2  &0.24e-6 & &4 &4 & 1.93 \\ 
                        & implicit lu &0.42e-6 & &4 &4 & 0.22 \\  \\
Schubert Broyden n=50   & mod.huang1  &0.30e-6 & &9 &9 & 4.73 \\
  $10x_0$               & mod.huang2  &0.83e-6 & &8 &8 & 3.90 \\ 
                        & implicit lu &0.48e-6 & &7 &7 & 0.44 \\  \\
Schubert Broyden n=50   & mod.huang1  &0.83e-6 & &11 &11 & 5.99 \\
  $100x_0$              & mod.huang2  &0.42e-6 & &12 &12 & 5.93 \\ 
                        & implicit lu &0.72e-6 & &11 &11 & 0.66 \\  \\
Schubert Broyden n=100  & mod.huang1  &0.24e-6 & &4 &4 & 16.32 \\
  $x_0$                 & mod.huang2  &0.24e-6 & &4 &4 & 15.49 \\ 
                        & implicit lu &0.42e-6 & &4 &4 & 1.26 \\  \\
Schubert Broyden n=100  & mod.huang1  &0.30e-6 & &9 &9 & 38.56 \\
  $10x_0$               & mod.huang2  &0.83e-6 & &8 &8 & 31.85 \\ 
                        & implicit lu &0.48e-6 & &7 &7 & 2.25 \\  \\
Schubert Broyden n=100  & mod.huang1  &0.83e-6 & &11 &11 & 47.62 \\
  $100x_0$              & mod.huang2  &0.42e-6 & &12 &12 & 46.41 \\ 
                        & implicit lu &0.72e-6 & &11 &11 & 3.57 \\  \\
\end{tabular}}
 
\newpage

{\scriptsize
\begin{tabular}{llccccc}  
\multicolumn{2}{c}{ }&\multicolumn{5}{c}{single precision } \\  \\
function & method &$\|F\|_\infty$& &$it_{ \|F\|_\infty}$&$it$&$time$\\  \\
Rosenbrock n=2 & mod.huang1  &0.0 & &1 &1 & 0.0 \\
  $x_0$        & mod.huang2  &0.0 & &1 &1 & 0.0 \\ 
               & implicit lu &0.0 & &1 &1 & 0.0 \\  \\
Rosenbrock n=2 & mod.huang1  &0.12e-6 & &1 &1 & 0.0 \\
$ 1.1x_0$      & mod.huang2  &0.12e-6 & &1 &1 & 0.0 \\ 
               & implicit lu &0.12e-6 & &1 &1 & 0.0 \\  \\
Rosenbrock n=2 & mod.huang1  &0.0 & &2 &2 & 0.0 \\
  $10x_0$      & mod.huang2  &0.0 & &2 &2 & 0.0 \\ 
               & implicit lu &0.0 & &1 &1 & 0.0 \\  \\
Rosenbrock n=2 & mod.huang1  &0.0 & &2 &2 & 0.0 \\
  $100x_0$     & mod.huang2  &0.0 & &2 &2 & 0.0 \\ 
               & implicit lu &0.0 & &1 &1 & 0.0 \\  \\
Extended Rosenbrock n=10 & mod.huang1  &0.0 & &1 &1 & 0.0 \\
  $x_0$                  & mod.huang2  &0.0 & &1 &1 & 0.0 \\ 
                         & implicit lu &0.0 & &1 &1 & 0.0 \\  \\
Extended Rosenbrock n=10 & mod.huang1  &0.12e-6 & &1 &1 & 0.0 \\
$ 1.1x_0$                & mod.huang2  &0.12e-6 & &1 &1 & 0.0 \\ 
                         & implicit lu &0.12e-6 & &1 &1 & 0.0 \\  \\
Extended Rosenbrock n=10 & mod.huang1  &0.0 & &2 &2 & 0.0 \\
  $10x_0$                & mod.huang2  &0.0 & &2 &2 & 0.0 \\ 
                         & implicit lu &0.0 & &1 &1 & 0.0 \\  \\
Extended Rosenbrock n=10 & mod.huang1  &0.0 & &2 &2 & 0.0 \\
  $100x_0$               & mod.huang2  &0.0 & &2 &2 & 0.0 \\ 
                         & implicit lu &0.0 & &1 &1 & 0.0 \\  \\
Extended Rosenbrock n=100 & mod.huang1  &0.0 & &1 &1 & 4.01 \\
  $x_0$                   & mod.huang2  &0.0 & &1 &1 & 2.96 \\ 
                          & implicit lu &0.0 & &1 &1 & 0.33 \\  \\
Extended Rosenbrock n=100 & mod.huang1  &0.12e-6 & &1 &1 & 4.01 \\
$ 1.1x_0$                 & mod.huang2  &0.12e-6 & &1 &1 & 2.96 \\ 
                          & implicit lu &0.12e-6 & &1 &1 & 0.33 \\  \\
Extended Rosenbrock n=100 & mod.huang1  &0.0 & &2 &2 & 8.02 \\
  $10x_0$                 & mod.huang2  &0.0 & &2 &2 & 5.94 \\ 
                          & implicit lu &0.0 & &1 &1 & 0.33 \\  \\
Extended Rosenbrock n=100 & mod.huang1  &0.0 & &2 &2 & 8.01 \\
  $100x_0$                & mod.huang2  &0.0 & &2 &2 & 5.93 \\ 
                          & implicit lu &0.0 & &1 &1 & 0.33 \\  \\
Powell singular     n=4   & mod.huang1  &0.64e-7  & &19 &19 & 0.0 \\
  $x_0$                   & mod.huang2  &0.81e-10 & &19 &19 & 0.0 \\ 
                          & implicit lu &0.38e+2 &(div) &1 &100 & 0.11 \\  
                   & imp.lu line search &0.40e-6 &  &86 &86 & 0.06 \\  \\
Powell singular     n=4   & mod.huang1  &0.44e-7 & &25 &25 & 0.0 \\
$ 1.1x_0$                 & mod.huang2  &0.56e-10 & &25 &25 & 0.0 \\ 
                          & implicit lu &0.42e+2 &(div) &1 &100 & 0.06 \\  
                   & imp.lu line search &0.61e-6 &  &60 &60 & 0.05 \\  \\
Powell singular     n=4   & mod.huang1  &0.34e-7 & &20 &20 & 0.0 \\
$ 10x_0$                  & mod.huang2  &0.43e-10 & &20 &20 & 0.0 \\ 
                          & implicit lu &0.38e+3 &(div) &1 &100 & 0.06 \\  
                   & imp.lu line search &0.69e-6 &  &69 &69 & 1.1 \\  \\
Powell singular     n=4   & mod.huang1  &0.70e-7 & &30 &30 & 0.0 \\
$ 100x_0$                 & mod.huang2  &0.88e-10 & &30 &30 & 0.0 \\ 
                          & implicit lu &0.93e+3 &(div) &9 &100 & 0.05 \\  
                   & imp.lu line search &0.69e-6 &  &73 &73 & 0.05 \\  
\end{tabular}}                                             

\newpage

{\scriptsize
\begin{tabular}{llccccc}  
\multicolumn{2}{c}{ }&\multicolumn{5}{c}{double precision } \\  \\
function & method &$\|F\|_\infty$& &$it_{ \|F\|_\infty}$&$it$&$time$\\  \\

Brown almost linear n=4 & mod.huang1  &0.22d-15 & &5 &5 & 0.0 \\
  $x_0$                 & mod.huang2  &0.44d-15 & &5 &5 & 0.0 \\ 
                        & implicit lu &0.11d-15 & &6 &6 & 0.0 \\  \\
Brown almost linear n=4 & mod.huang1  &0.0      & &5 &5 & 0.0 \\
  $1.1x_0$              & mod.huang2  &0.0      & &5 &5 & 0.0 \\ 
                        & implicit lu &0.0      & &6 &6 & 0.0 \\  \\
Brown almost linear n=4 & mod.huang1  &0.67d-15 & &8 &8 & 0.0 \\
  $10x_0$               & mod.huang2  &0.22d-15 & &5 &5 & 0.0 \\ 
                        & implicit lu &0.10d+1  & &1 &2 & 0.0 \\  \\
Brown almost linear n=4 & mod.huang1  &0.0     & &16 &16 & 0.0 \\
    $100x_0$     & m.hua1 line search &0.78d-15& &11 &11 & 0.0 \\  
                        & mod.huang2  &0.11d-15& &16 &16 & 0.0 \\ 
                 & m.hua2 line search &0.33d-15& &11 &11 & 0.0 \\ 
                 & implicit lu &0.44d-15 & &18 &18 & 0.0 \\  \\
Brown almost linear n=20& mod.huang1  &0.56d-14 &(o) &5 &9 & 0.44 \\
    $x_0$        & m.hua1 line search &0.44d-15 & &21  &21  & 0.93 \\  
                        & mod.huang2  &0.24d-14 &(o) &15 &22 & 1.05 \\ 
                 & m.hua2 line search &0.89d-15 & &22  &22  & 1.05 \\ 
                        & implicit lu &0.29d-14 &(o) &12 &17 & 0.22 \\  
                 & imp.lu line search &0.89d-15 & &14 &14 & 0.16 \\ \\ 
Brown almost linear n=20& mod.huang1  &0.29d-14 &(o) &16 &21 & 0.88 \\
    $1.1x_0$     & m.hua1 line search &0.89d-15 & &22  &22  & 0.99 \\  
                        & mod.huang2  &0.40d-14 &(o) &9 &14 & 0.66 \\ 
                 & m.hua2 line search &0.78d-15 & &12  &12  & 0.60 \\ 
                        & implicit lu &0.32d-14 &(o) &8 &12 & 0.11 \\ 
                 & imp.lu line search &0.44d-15 & &16 &16 & 0.17 \\  \\
Schubert Broyden n=10   & mod.huang1  &0.99d-15 & &5 &5 & 0.0 \\
  $x_0$                 & mod.huang2  &0.99d-15 & &5 &5 & 0.0 \\ 
                        & implicit lu &0.88d-15 & &5 &5 & 0.0 \\  \\
Schubert Broyden n=10   & mod.huang1  &0.89d-15 & &9 &9 & 0.5 \\
  $10x_0$               & mod.huang2  &0.89d-15 & &9 &9 & 0.5 \\ 
                        & implicit lu &0.89d-15 & &8 &8 & 0.0 \\  \\
Schubert Broyden n=50   & mod.huang1  &0.78d-15 & &5 &5 & 2.80 \\
  $x_0$                 & mod.huang2  &0.78d-15 & &5 &5 & 2.63 \\ 
                        & implicit lu &0.89d-15 & &5 &5 & 0.33 \\  \\
Schubert Broyden n=50   & mod.huang1  &0.67d-15 & &10 &10 & 5.55 \\
  $10x_0$               & mod.huang2  &0.99d-15 & &9 &9 & 4.78 \\ 
                        & implicit lu &0.67d-15 & &9 &9 & 0.60 \\  \\
Schubert Broyden n=50   & mod.huang1  &0.78d-15 & &13 &13 & 7.25 \\
  $100x_0$              & mod.huang2  &0.78d-15 & &13 &13 & 6.92 \\ 
                        & implicit lu &0.89d-15 & &12 &12 & 0.77 \\  \\
Schubert Broyden n=100  & mod.huang1  &0.89d-15 & &5 &5 & 21.37 \\
  $x_0$                 & mod.huang2  &0.89d-15 & &5 &5 & 20.71 \\ 
                        & implicit lu &0.89d-15 & &5 &5 & 1.70 \\  \\
Schubert Broyden n=100  & mod.huang1  &0.89d-15 & &10 &10 & 42.74 \\
  $10x_0$               & mod.huang2  &0.89d-15 & &10 &10 & 41.57 \\ 
                        & implicit lu &0.89d-15 & &9 &9 & 3.08 \\  \\
Schubert Broyden n=100  & mod.huang1  &0.89d-15 & &13 &13 & 55.64 \\
  $100x_0$              & mod.huang2  &0.78d-15 & &13 &13 & 54.10 \\ 
                        & implicit lu &0.89d-15 & &12 &12 & 4.12 \\  \\
\end{tabular}}

\newpage

{\scriptsize
\begin{tabular}{llccccc}  
\multicolumn{2}{c}{ }&\multicolumn{5}{c}{double precision } \\  \\
function & method &$\|F\|_\infty$& &$it_{ \|F\|_\infty}$&$it$&$time$\\  \\
Rosenbrock n=2 & mod.huang1  &0.22d-15 & &1 &1 & 0.0\\
     $x_0$     & mod.huang2  &0.22d-15 & &1 &1 & 0.0 \\ 
               & implicit lu &0.22d-15 & &1 &1 & 0.0 \\  \\
Rosenbrock n=2 & mod.huang1  &0.22d-15 & &1 &1 & 0.0\\
 $1.1x_0$      & mod.huang2  &0.22d-15 & &1 &1 & 0.0 \\ 
               & implicit lu &0.22d-15 & &1 &1 & 0.0 \\  \\
Rosenbrock n=2 & mod.huang1  &0.0 & &1 &1 & 0.0\\
$10x_0$        & mod.huang2  &0.0 & &1 &1 & 0.0 \\ 
               & implicit lu &0.0 & &1 &1 & 0.0 \\  \\
Rosenbrock n=2 & mod.huang1  &0.0 & &1 &1 & 0.0\\
  $100x_0$     & mod.huang2  &0.0 & &1 &1 & 0.0 \\ 
               & implicit lu &0.0 & &1 &1 & 0.0 \\  \\
Extended Rosenbrock n=10 & mod.huang1  &0.22d-15 & &1 &1 & 0.0 \\
  $x_0$                  & mod.huang2  &0.22d-15 & &1 &1 & 0.0 \\ 
                         & implicit lu &0.22d-15 & &1 &1 & 0.0 \\  \\
Extended Rosenbrock n=10 & mod.huang1  &0.22d-15 & &1 &1 & 0.0 \\
$ 1.1 x_0$               & mod.huang2  &0.22d-15 & &1 &1 & 0.0 \\ 
                         & implicit lu &0.22d-15 & &1 &1 & 0.0 \\  \\
Extended Rosenbrock n=10 & mod.huang1  &0.0 & &1 &1 & 0.0 \\
  $10x_0$                & mod.huang2  &0.0 & &1 &1 & 0.0 \\ 
                         & implicit lu &0.0 & &1 &1 & 0.0 \\  \\
Extended Rosenbrock n=10 & mod.huang1  &0.0 & &1 &1 & 0.0 \\
  $100x_0$               & mod.huang2  &0.0 & &1 &1 & 0.0 \\ 
                         & implicit lu &0.0 & &1 &1 & 0.0 \\  \\
Extended Rosenbrock n=100 & mod.huang1  &0.22d-15 & &1 &1 & 4.23 \\
  $x_0$                   & mod.huang2  &0.22d-15 & &1 &1 & 3.18 \\ 
                          & implicit lu &0.22d-15 & &1 &1 & 0.33 \\  \\
Extended Rosenbrock n=100 & mod.huang1  &0.22d-15 & &1 &1 & 4.23 \\
$ 1.1 x_0$                & mod.huang2  &0.22d-15 & &1 &1 & 3.13 \\ 
                          & implicit lu &0.22d-15 & &1 &1 & 0.33 \\  \\
Extended Rosenbrock n=100 & mod.huang1  &0.0 & &1 &1 & 4.29 \\
  $10x_0$                 & mod.huang2  &0.0 & &1 &1 & 3.13 \\ 
                          & implicit lu &0.0 & &1 &1 & 0.28 \\  \\
Extended Rosenbrock n=100 & mod.huang1  &0.0 & &1 &1 & 4.23 \\
  $100x_0$                & mod.huang2  &0.0 & &1 &1 & 3.14 \\ 
                          & implicit lu &0.0 & &1 &1 & 0.38 \\  \\
Powell singular     n=4   & mod.huang1  &0.62d-16  & &44 &44 & 0.0 \\
  $x_0$                   & mod.huang2  &0.78d-19 & &44 &44 & 0.0 \\ 
                          & implicit lu &0.38d+2 &(div) &1 &100 & 0.0 \\  
                   & imp.lu line search &0.56d-15 &  &186 &186 & 0.11 \\  \\
Powell singular     n=4   & mod.huang1  &0.33d-16  & &45 &45 & 0.0 \\
  $1.1x_0$                & mod.huang2  &0.41d-19 & &45 &45 & 0.0 \\ 
                          & implicit lu &0.42d+2 &(div) &1 &100 & 0.05 \\  
                   & imp.lu line search &0.85d-15 &  &126 &126 & 0.06 \\  \\
Powell singular     n=4   & mod.huang1  &0.43d-16  & &50 &50 & 0.0 \\
  $10x_0$                 & mod.huang2  &0.54d-19 & &50 &50 & 0.0 \\ 
                          & implicit lu &0.38d+3 &(div) &1 &100 & 0.05 \\  
                   & imp.lu line search &0.56d-15 &  &135 &135 & 0.11 \\  \\
Powell singular     n=4   & mod.huang1  &0.70d-16  & &55 &55 & 0.0 \\
  $100x_0$                & mod.huang2  &0.84d-19 & &55 &55 & 0.0 \\ 
                          & implicit lu &0.93d+4 &(div) &9 &100 & 0.05 \\  
                   & imp.lu line search &0.39d-15 &  &140 &140 & 0.05 \\  \\
\end{tabular}}                                             
\end{document}